  \documentclass[12pt]{article}
  \usepackage{amsfonts}

   \def\R{\mathbb{R}}
   \def\N{\mathbb{N}}  
     
   \def\Q{\mathbb{Q}}  

   \def\e{{\varepsilon}}
   \def\D{{\nabla}}
   \def\vi{{\varphi}}

   \def\cB{{\cal B}}
   \def\cC{{\cal C}}
   \def\cD{{\cal D}} 
   \def\cE{{\cal E}}

   \def\cU{{\cal U}}

   \def\dist{\mathop{\rm dist}\nolimits}
   \def\lo{\mathop{\longrightarrow}}
   
   \def\loc{\mathop{\rm loc}\nolimits}

   \def\ccup{\mathop{\cup}}
   \def\const{\mathop{\rm const}\nolimits}  
   
   \def\qed{\hfill{\em q.e.d.}\\ \vspace{1mm}}
   \def\proof{\noindent{\underline {\sf Proof} \hspace{2mm}}}

\newtheorem{df}{Definition}[section]
\newtheorem{prop}[df]{Proposition}
\newtheorem{lemma}[df]{Lemma}
\newtheorem{teo}[df]{Theorem}
\newtheorem{rem}[df]{Remark}

\newtheorem{cor}[df]{Corollary}

 \newcommand{\sezione}[1]{\section{#1}\setcounter{equation}{0}}

\begin{document}


   \title{Positive solutions for a nonlinear elliptic problem with
     strong lack of compactness\thanks{Work supported by the Italian
national research  project ``Metodi 
variazionali e topologici nello studio di fenomeni non lineari".}
\author{Riccardo Molle}}

\date{}
 
 \maketitle


\vspace{-3mm}

 \begin{center}

{\em
Dipartimento di Matematica\\
Universit\`a di Roma ``Tor Vergata''\\
Via della Ricerca Scientifica n. 1\\
00133 Roma, Italy\\
  E-mail: {\sf  molle@mat.uniroma2.it}
}

\end{center}

\vspace{5mm}
 

 \begin{abstract}

This paper deals with the lack of compactness in nonlinear elliptic
problems $(P)$. 
In particular, a domain $\Omega$ is provided where not converging
Palais-Smale sequences exist at every energy level.
Nevertheless, it is proved that problem $(P)$ has infinitely many
solutions on $\Omega$.

 \end{abstract}

{\small


\noindent{\bf A.M.S. subject classification 2000:} 35J20; 35J65.
}


\sezione{Introduction, examples and statement of the result}


In this paper we are concerned with the failure of compactness related
to problem
\[
(P) \hspace{1cm}
\left\{
        \begin{array}{ll}
         -\Delta u+u=|u|^{p-2}u          & \mbox{ in }\Omega,\\
         u>0                            & \mbox{ in }\Omega,\\
         u \in H_0^1(\Omega) 
        \end{array}
\right.
\]
when $\Omega\subset \R^N$, $N\ge 3$, is an unbounded domain,
having unbounded boundary, and $2<p<2^*=2N/(N-2)$. 

There are many studies related to problem $(P)$, in recent years,
because of its connection with questions arising from applied
sciences.
A typical approach is the variational one: the  solutions of Problem $(P)$  can be found in correspondence to
positive functions that are critical points of the ``energy''
functional $\cE:M\to\R$ given by 
\[
\cE(u)=\int_\Omega(|\D u|^2+u^2)dx,
\]
where $M$ is the manifold
\[
M=\{u\in H_0^1(\Omega)\ :\ |u|_{L^p(\Omega)}=1\}.
\]

When $\Omega$ is bounded, classical variational methods can be applied,
taking advantage of the compactness of the embedding
$H^1_0(\Omega)\hookrightarrow L^p(\Omega)$, so problem $(P)$ always
has a solution, corresponding to $\min_M\cE$, and the multiplicity of
solutions is related to the shape of $\Omega$ (see, for example,
\cite{BCP,We} and references therein).

The key compactness condition is the {\bf Palais-Smale condition}:
{\em $\cE$ satisfies the Palais-Smale condition  if every sequence
  $(u_n)_n$ in $M$ such that $\lim\limits_{n\to+\infty}\cE(u_n)=c$,
  $c\in\R$, $\lim\limits_{n\to +\infty}\D\cE (u_n)=0$ (a Palais Smale
  sequence), has a converging subsequence.}  

If $\Omega$ is unbounded the situation is different from the bounded case. 
Of course, when the embedding $H^1_0(\Omega)\hookrightarrow
L^p(\Omega)$ is compact, it is again possible to apply the variational
techniques used in the bounded case.
This happens, for example, if $\Omega$ is ``{\em thin}'' at infinity,
in the sense that
\[
\lim_{R\to+\infty}\sup\{\mu (B(x,1)\cap\Omega)\ :\ x\in\R^N,\ |x|=R\}=0,
\]
where $\mu$ denotes the Lebesgue measure and $B(x,1)$ is the unitary
ball with centre in $x$. 
So, let us consider the case when $\Omega$ is not thin at infinity.
If $\Omega=\R^N$, $(P)$ has a solution $\omega$, unique up to
translations, corresponding to  
\begin{equation}
\label{em}
m:=\min\{\|u\|^2_{H^1(\R^N)}\ :\ u\in H^1(\R^N),\ |u|_{L^p(\R^N)}=1\}.
\end{equation}
When $\Omega$ looks like $\R^N$ at infinity, in the sense that
$\R^N\setminus\Omega$ is bounded (i.e. $\Omega$ is an {\em exterior domain}),
then the analysis of the compactness failure shows that a Palais-Smale
sequence of $\cE$ differs from its weak limit by ``waves'' $\omega$
that go to infinity.
This behaviour of the Palais-Smale sequences ensures that a local
compactness condition holds.
Taking into account this property, existence and multiplicity results
for problem $(P)$ in exterior domains are proved (see \cite{BC,BL,MP}
and references therein).

When not only $\Omega$, but also $\R^N\setminus\Omega$ is unbounded,
the existence of solutions of $(P)$ depends on the shape of
$\Omega$. For example in \cite{EL} it is proved that if  there exists $\bar
x\in\R^N$ such that $(\nu (x),\bar x)\ge 0\ \forall x\in\partial\Omega$, $(\nu
(x), \bar x)\not\equiv 0$, then problem $(P)$ has no solution, while
in \cite{E} it is shown that solutions for $(P)$ exist in
``strip-like'' domains (see also the monograph \cite{Wa}, and
references therein, for other existence results in domains with
unbounded boundary).

In \cite{M} a general condition on the shape of $\Omega$ at infinity
is stated such that, even if both $\Omega$ and $\R^N\setminus\Omega$
are unbounded, a local compactness condition still holds.
The condition stated in \cite{M} requires that $\Omega$ enlarges, at
infinity, and its boundary flattens, or shrinks, at infinity. 
Namely, the following conditions are required:
\[
\left\{ \begin{array}{ll}
\vspace{2mm}
{\displaystyle \lim_{R\to +\infty} \inf\{r(x)\ :\ x\in\Omega,\
|x|=R\}=+\infty,}&\hspace{2cm}(C_1)\\
{\displaystyle \lim_{R\to +\infty}\sup\{ h(y)\ :\
y\in\partial\Omega,\ |y|=R\}=0} &\hspace{2cm}(C_2)
\end{array}\right.
\]
where
\[
r(x)=\sup\{\rho>0\ : \ \exists\bar x\in\Omega\ \mbox{ such that } x\in
B(\bar x,\rho)\mbox{ and } B(\bar x,\rho)\subset\Omega\},
\]
\begin{eqnarray*}
h(y)=\sup\{\dist(z,T_{\partial\Omega,y}\cap B(y,1))&:& z\mbox{ is in
the connected component}\\
& & \mbox{ of }\partial\Omega\cap B(y,1)\mbox{ containing } y\}
\end{eqnarray*}
(here $T_{\partial\Omega,y}$ is the hyperplane tangent to $\partial\Omega$
in $y$).
It is worth to remark that conditions $(C_1)$ and $(C_2)$, hence local
compactness, is not sufficient to
guarantee the existence of solutions of $(P)$; for example, an half-space $\Pi$
verifies $(C_1)$ and $(C_2)$, but $(P)$ has no solution on $\Pi$, by
the Esteban-Lions result.
What happens in half-spaces is that, actually, no Palais-Smale
sequences for $\cE$ exist at the compactness levels.

Domains $\Omega$ with unbounded boundary and periodic in some
directions (for example the exterior of a cylinder) are considered in
\cite{CMP} and a result of multiplicity of solutions for $(P)$ is
proved.
Let us remark that every solution $\bar u$, found in \cite{CMP},
provides a non-compact family of solutions of $(P)$ (hence a
non-compact Palais-Smale sequence at the level $\cE\left({\bar u\over
  |\bar u|_{L^p(\Omega)}}\right)$) by moving $\bar u$ according to the
periodic structure of $\Omega$.

Now, one can look for a domain $\Omega$ with no periodicity assumption
and such that not converging Palais-Smale sequences for $\cE$ exist at
every energy level.
It is worth to observe that in a domain with this property no local
compactness condition can hold. 

Looking in finite dimension, a function $f\in
\cC^1([0,+\infty),(0,+\infty))$ such that
\begin{equation}
\label{t1}
f(i)=q_i,\quad f'(i)=0\qquad\forall i\in\N,\quad\mbox{ where
}(q_i)_i=\Q^+,
\end{equation}
has such a strong failure of compactness.
Indeed for every $s\in [0,+\infty)$ we can
    consider $(q_{i_j})_j$ in $\Q^+$ such that $\lim\limits_{j\to
      +\infty}f(i_j)=q_{i_j}=s$ and so, taking into account
    (\ref{t1}), we have that $(i_j)$ is a not converging Palais-Smale
    sequence for $f$ at the level $s$.
Observe that every $q\in\Q$ is a critical level for $f$, i.e. the
critical levels are a dense subset of the range $(0,+\infty)$ (no more can be
required by the Sard's Theorem).
It is also possible to consider a function $g\in
\cC^1(\R^2,(0,+\infty))$ such that
\begin{equation}
\label{t2}
g(i,0)=q_i,\quad \D g(i,0)=\left({1\over i}\, ,\, {1\over
  i}\right),\quad{\partial g\over\partial y}\, (x,y)>0\quad \forall
  (x,y)\in\R^2.
\end{equation}
In this second example we see that there is the same strong failure of
the compactness condition of the previous example and $g$ has no
critical point.

\vspace{1mm}

Going back to the functional $\cE$, define
\begin{equation}
\label{S}
S_i=\left\{(x_1,\ldots,x_N)\in\R^N\ :\ -\mbox{${1\over
    2}q_i$}<x_N<\mbox{${1\over 2}q_i$}\right\}, 
\end{equation}
 ($(q_i)_i=\Q^+$) and set
\begin{eqnarray*}
\Omega & = &\left\{x=(x_1,\ldots,x_N)\in\R^N\ :\
\mbox{$\sum_{j=1}^{N-1} x_j^2<1$}\right\}\cup\\
& & S_1
\cup\left[\ccup_{i=2}^{\infty}[S_i+(0,\ldots,0,i+\mbox{$\sum_{j=1}^{i-1}q_j$}+
\mbox{${1\over 2}q_i$})]\right].
\end{eqnarray*}
The domain $\Omega$ is introduced in \cite{MM}, where it is seen that
the following properties hold.
\begin{prop}
\label{MM}
It holds:

$a)\ \ \{\cE(u)\ :\ u\in M\}=(m,+\infty)$;

$b)\  $ for every $c\in[m,+\infty)$ there exists a Palais-Smale
  sequence $(u_n)_n$ for $\cE$ at the level $c$, i.e.
\[
\lim_{n\to +\infty} \cE(u_n)=c,\qquad \lim_{n\to +\infty}\D \cE(u_n)=0,
\]
that is not relatively compact.
\end{prop}

For the reader's convenience, a proof of Proposition \ref{MM} is
contained in Section 2.
Here we stress the fact that not converging Palais-Smale sequence for
$\cE$ are related either to translations in the $x_i$-directions, for
$i=1,\ldots,N-1$,  either to translations in the $x_N$-direction. 

\vspace{1mm}

In this paper we prove that, in spite of the strong failure in
compactness described above and of the non periodic structure of
$\Omega$, problem $(P)$ on $\Omega$ has infinitely many solutions. 

\begin{teo}
\label{T}
Problem $(P)$ has a sequence of solutions $(v_n)_n$ such that
\[
\lim_{n\to +\infty}\cE\left({v_n\over |v_n|_{L^p(\Omega)}}\right)=m.
\]
\end{teo}

The solutions given by Theorem \ref{T} correspond to local minima,
that are ``localized'' in large strips.
We reach the local minima by following maximal slope curves.
In order to construct such curves, a basic tool is the analysis of the
Palais-Smale sequences (\S 3).
Indeed, we show that Palais-Smale sequences are compact on a maximal
slope curve, near the local minima.

\vspace{1mm}

Finally, let us remark that if we consider the domain 
\[
\cD=\Omega\cup\{x=(x_1,\ldots,x_N)\in\R^N\ :\ x_1<0\}
\]
then Proposition \ref{MM} still holds, on $\cD$, but problem $(P)$ has
no solution, because $\cD$ satisfies the Esteban-Lions condition. 

\vspace{1mm}

The paper is organized as follows: Section 2 contains the proof of
Proposition \ref{MM}, the variational
setting and some useful tools; Section 3 is devoted to
the analysis of the Palais-Smale sequences and in Section 4 Theorem
\ref{T} is proved.


\sezione{Preliminary results, the variational setting and useful tools}


In this paper the following notations are used:
{\small
\begin{itemize}
\item
if $\cD\subseteq \R^N$ and $u\in H_0^1(\cD)$, we
denote also by $u$ its extension to $\R^N$ obtained by setting $u\equiv
0$ outside $\cD$.
\item
$L^p(\cD)$, $1\le p < +\infty$, $\cD\subseteq\R^N$, denotes the usual 
  Lebesgue space, endowed with the norm $|u|_{L^p(\cD)}=(\int_\cD
  |u|^p)^{1\over p}$; if $\cD=\R^N$, we simply write $|u|_{L^p(\R^N)}=|u|_p$.
\item
$H^1_0(\cD)$, $\cD\subseteq \R^N$
  denotes the Sobolev space obtained as the closure of
  $\cC_0^\infty(\cD)$ with respect to the
  norm
\[
\|u\|_{H^1(\cD)}=\left[\int_\cD(|\D u|^2+u^2)dx\right]^{1\over 2};
\]
if $\cD=\R^N$, we simply write $\|u\|_{H^1(\R^N)}=\|u\|$.
\item
The generic point
$x=(x_1,\ldots,x_{N-1},x_N)\in\R^N$
is denoted by $(x',x_N)$, where 
$x'=(x_1,\ldots,x_{N-1})\in \R^{N-1}$ and $x_N\in\R$; we put also
$|x'|=\left(\sum\limits_{j=1}^{N-1} x_j^2\right)^{1\over 2}$.
\item
We set the points $Q_1^-=(0,\ldots,0,-{q_1\over2})$, $Q_1^+=(0,\ldots,0,{q_1\over2})$
 and, for $i=2,3,\ldots$,
 $Q^-_i=(0,\ldots,0,i+\mbox{$\sum_{j=1}^{i-1}q_j$})$,\qquad
$Q^+_i=(0,\ldots,0,i+\mbox{$\sum_{j=1}^{i}q_j$})$.
\end{itemize}
}

\vspace{2mm}

\noindent{\underline {\sf Proof of Proposition \ref{MM}.} \hspace{2mm}}
First, we prove that
\begin{equation}
\label{eem}
\inf_M \cE=m.
\end{equation}
Let us recall that $m$ in (\ref{em}) is achieved by a positive
function $\omega$, radially symmetric and decreasing when the radial
co-ordinate increases, that is unique modulo translations and satisfies
\begin{equation}
\label{as1}
\lim_{|x|\to+\infty}|\omega(x)||x|^{\frac{N-1}{2}}e^{|x|}=d>0,
\end{equation}
\begin{equation}
\label{as2}
\lim_{|x|\to+\infty}|\D \omega(x)||x|^{\frac{N-1}{2}}e^{|x|}=d,
\end{equation}
for a suitable positive constant $d$ (see \cite{BL} and references therein). 

For $R>0$, let $\vi_R:\R^N\to \R$ be a cut-off function defined by
$\vi_R(x)=\vi\left({|x|\over R}\right)$, where $\vi\in
\cC^\infty(\R^+,[0,1])$ is a non-increasing function such that
$\vi(t)=1$ if $t\in \left[0,{1\over 4}\right]$, $\vi(t)=0$ if
$t\in\left[{1\over 2},+\infty\right)$, and define
\[
u_i(x)=\frac
{\vi_{q_i}\left(x-{Q_i^-+Q_i^+\over
    2}\right)\omega\left(x-{Q_i^-+Q_i^+\over 2}\right)}
{\left|\vi_{q_i}\left(x-{Q_i^-+Q_i^+\over
    2}\right)\omega\left(x-{Q_i^-+Q_i^+\over 2}\right)\right|_p}.
\]
It is clear that $u_i\in M$, $\forall i\in \N$, and, taking into
account (\ref{as1}) and (\ref{as2}), we have
\begin{equation}
\label{mm2}
\lim_{q_i\to +\infty} \left|\vi_{q_i}\left(x-\mbox{${Q_i^-+Q_i^+\over
    2}$}\right)\omega\left(x-\mbox{${Q_i^-+Q_i^+\over 2}$}\right)-
\omega\left(x-\mbox{${Q_i^-+Q_i^+\over 2}$}\right) \right|_p=0,
\end{equation}
\begin{equation}
\label{mm3}
\lim_{q_i\to +\infty} \left\|\vi_{q_i}\left(x-\mbox{${Q_i^-+Q_i^+\over
    2}$}\right)\omega\left(x-\mbox{${Q_i^-+Q_i^+\over 2}$}\right)
-\omega\left(x-\mbox{${Q_i^-+Q_i^+\over 2}$}\right)\right\|=0.
\end{equation}
From (\ref{mm2}) and (\ref{mm3}) we infer $\lim\limits_{q_i\to
  +\infty}\left\|u_i(x)-\omega\left(x-\mbox{${Q_i^-+Q_i^+\over
    2}$}\right)\right\|=0$, i.e.
\begin{equation}
\label{mm1}
\lim_{q_i\to +\infty} \cE(u_i)=m,
\end{equation}
that implies $\inf\limits_M \cE\le m$. On the other hand  $\inf\limits_M
\cE\ge m$ directly follows from (\ref{em}), so (\ref{eem}) is proved.

\vspace{1mm}

We claim that the infimum in (\ref{eem}) is not achieved.
Assume, by contradiction, that a minimum point $u^*\in M$ exists.
By the uniqueness of the minimizers of (\ref{em}), there exists
$y^*\in\R^N$ such that $u^*(x)=\omega(x-y^*)$, $\forall x\in\R^N$.
This is not possible because $u^*\equiv 0$ on $\R^N\setminus\Omega$,
while $\omega>0$ on $\R^N$, so we have proved the claim.

\vspace{1mm}

Recall, now, that for every $i\in\N$ there exists a critical point
$z_i$ for $\cE$ on $\{u\in H^1_0(S_i)\ :\ |u|_p=1\}$ ($S_i$ is
introduced in (\ref{S})), corresponding to
the ``minimal'' solution of the problem $(P)$ on $S_i$ (see \cite{E}). 
Namely, $z_i$ satisfies
\begin{equation}
\label{q1}
\Theta(q_i):=\|z_i\|^2=\min\{\|u\|^2\ :\ u\in H^1_0(S_i),\  |u|_p=1\}.
\end{equation}

It is easily seen that the above introduced $\Theta:\Q^+\to \R$ is a
continuous monotone decreasing function and we claim that
\begin{equation}
\label{q00}
 \lim_{q_i\to+\infty}\Theta(q_i)  =m,
\end{equation}
\begin{equation}
\label{q0}
\lim_{q_i\to 0}\Theta(q_i) =+\infty.
\end{equation}
Taking into account (\ref{q1}), in order to show (\ref{q00}) we can
use the test functions ${\vi_{q_i}\omega\over |\vi_{q_i}\omega|_p}$ and
  proceed with the same argument developed to prove (\ref{eem}).

To prove (\ref{q0}) observe that, by the Poincar\'e's inequality on a
strip, there exists a constant $\bar k>0$ such that
\begin{equation}
\label{q3}
|u|_{p}\leq\bar k\, |\D u|_{L^2(S_{1})}\qquad\forall u\in H^1_0( S_1)
\end{equation}
(see Theorem 6.30 in \cite{A}, for example). To simplify the notation,
let us assume $q_1=1$ and set $\hat z_i(x)=z_i(q_ix)$; applying
(\ref{q3}) to $\hat z_i$ and taking into account that $|z_i|_{p}=1$
$\forall i\in\N$,  we get
\[
 q_i^{-{(N-2)\over 2p}\left({2N\over N-2}-p\right)}\leq\bar k \, |\D
 z_i|_{L^2(S_{i})}<\bar k\, \|z_i\|,
\]
that implies (\ref{q0}).

Now we are in position to verify $(b)$.
If $c=\Theta (q_i)$, for some $i\in\N$, let us construct the
sequence $(z_{i,j})_j$  by
\[
z_{i,j}(x_1,\ldots,x_N)=z_i \left(x_1+j,x_2,\ldots,x_{N-1}, x_N-\mbox{${Q_i^-+Q_i^+\over 2}$}\right).
\]
Then, it is not difficult to see that
$(z_{i,j})_j$ is a Palais-Smale sequence for $\cE$ at the level $c$, and
it cannot have converging subsequence because it weakly converges to zero.

If $c\in[m,+\infty)\setminus\Theta(\Q^+)$, let $(q_{i_j})_j$ be a
  sequence in $\Q$ such that
  $c=\lim\limits_{j\to+\infty}\Theta(q_{i_j})$. Then $(z_{{i_j},j})_j$ is a
  Palais-Smale sequence  for $\cE$ at the level $c$ that has no converging
  sub\-sequences.

\qed

From the minimality of $m$ (see (\ref{em})) and $(a)$ of Proposition
\ref{MM} we get the following results, which we need in the sequel.

\begin{prop}
\label{C1}
Let $\mu>0$ and $\cD\subseteq \R^N$ be an open domain with piecewise
smooth boundary (if not empty).
If $\bar u$ is a nontrivial solution of 
\begin{equation}
\label{PmD}
P(\mu,\cD)  \hspace{1cm}
\left\{
        \begin{array}{l}
         -\Delta u+u=\mu |u|^{p-2}u          \quad \mbox{ in }\cD,\\
         u\in H^1_0(\cD),
        \end{array}
\right.
\end{equation}
then
\begin{equation}
\label{e2.7}
|\bar u|_{p}>\left({m\over \mu}\right)^{1\over p-2}\quad\mbox{ if
 }\cD\neq\R^N,\qquad|\bar u|_{p}\ge\left({m\over \mu}\right)^{1\over
 p-2}\quad\mbox{ if }\cD =\R^N.
\end{equation}
\end{prop}

\proof By (\ref{em}), 
\begin{equation}
\label{eC1}
m|\bar u|_{p}^2\leq \|\bar u\|^2,
\end{equation}
and, as we have seen in the proof of Proposition \ref{MM}, the equality in
(\ref{eC1}) can hold only if $\cD =\R^N$. 
Moreover, $\bar u$ being a solution of $P(\mu,\cD)$,
\[
\|\bar u\|^2=\mu |\bar u|^p_{p}.
\]
Thus
\[
|\bar u|_{p}^{p-2}\ge {m\over \mu},
\]
where   equality  holds only if $\cD =\R^N$,
so (\ref{e2.7}) follows.

\qed

\begin{cor}
\label{C2}
Let $\bar u$ be a critical point of $\cE$.
If $\cE(\bar u)\in (m,2^{1-{2\over p}}m]$ then $\bar u$ does not change
  sign.
\end{cor}

\proof Let $\bar u=\bar u^+-\bar u^-$,  with $\bar u^+\not\equiv 0$
and $\bar u^-\not\equiv 0$, and call $\cE(\bar u)=\mu$.
Then,  by $(a)$ in Proposition \ref{MM} and since $\bar u$ solves
$P(\mu,\Omega)$, we have
\begin{equation}
\label{eC2}
m|\bar u^\pm|_{p}^2 <\|\bar u^\pm\|^2=\mu |\bar
u^\pm|^p_{p}.
\end{equation}
From (\ref{eC2}) and $|\bar u|_p=1$ we infer
\[
1=|\bar u|_{p}^p=|\bar u^+|_{p}^p+|\bar u^-|_{p}^p
> 2\left({m\over \mu}\right)^{p\over p-2},
\]
that implies
\[
\mu >  2^{1-{2\over p}}m.
\]

\qed

\vspace{1mm}

In view of the axial symmetry of the problem and taking into account the
principle of symmetric criticality (see \cite{Pal}, or also \cite{W},
Theorem 1.28), to solve problem $(P)$ it is sufficient to consider the
subspace of the axially symmetric functions
\[
H^1_{0,r}(\Omega)=\{u\in H^1_0(\Omega)\ :\ u(x',x_N)=u(\bar x',x_N)\
\mbox{ if } |x'|=|\bar x'|\},
\]
the smooth manifold
\[
V=\{u\in H^1_{0,r}(\Omega)\ :\ |u|_p=1\}
\]
and look for critical points of the functional $E=\cE_{\mid V}$.

\begin{rem}
\label{R5}
{\em Let us observe that we cannot find critical points of $E$ simply by
minimization.
Indeed, the argument developed in the proof of $(a)$ of Proposition
\ref{MM} shows that
\begin{equation}
\label{eC5}
\inf_V E=m
\end{equation}
and that the infimum in (\ref{eC5}) is not achieved.
}\end{rem}

\vspace{2mm}

Now, we denote the $x_N$-axis by
\[
A=\{x=(x_1,\ldots,x_N)\in\R^N\ :\ x_1=\ldots=x_{N-1}=0\}.
\]
Then, we introduce a barycenter type function on $H^1_{0,r}(\R^N)\setminus\{0\}$: for
$u\in H^1_{0,r}(\R^N)\setminus\{0\}$ set 
\[
\tilde u(x)={1\over \mu(B(x,1))}\int_{B(x,1)} |u(y)|\,dy\qquad\forall x\in\R^N,
\]
\[
\hat u(x)=\left[\tilde u(x)-{1\over 2}\max_{\R^N} \tilde u(x)
  \right]^+
\qquad\forall x\in\R^N
\]
and define $\beta:H^1_{0,r}(\R^N)\setminus\{0\} \to A $ by
\begin{equation}
\label{e4.1}
\beta(u)={1\over |\hat u|^p_{p}}\int_{\R^N} (\hat u(x))^p x \,  dx.
\end{equation}

Let us  remark that $\beta$ is well defined for all $u\in H^1_{0,r}(\R^N)\setminus\{0\}$, because $\hat
u\not\equiv 0$ and has compact support, $\beta(u)\in A$ because $u$ is
axially symmetric and, moreover, $\beta$ is continuous and verifies:
\begin{equation}
\label{b1}
 \beta(u(x-y))=\beta(u(x))+y\qquad\forall u\in H^1_{0,r}(\R^N)\setminus\{0\},\ \forall y\in A,
\end{equation}
\begin{equation}
\label{b2}
\beta (\omega(x))=0.
\end{equation}


\sezione{Behaviour of Palais-Smale sequences of $E$}


In the next Proposition we study Palais-Smale sequences in the energy
range $[m,$ $2^{1-{2\over p}}m)$ (see Remark \ref{R} for the level $2^{1-{2\over p}}m$).
We show that or such a sequence is relatively compact or it
is composed by a unique ``wave'' that goes to infinity, in the
$x_N$-direction, and converges to a limit problem.

\begin{prop}
\label{PS}
Let $(u_n)_n$ be a Palais-Smale sequence for $E$ at level $c$.
If $c < 2^{1-{2\over p}}m$ then, up to a subsequence, one of the
following alternatives is possible:
\begin{itemize}
\item [$a)$]
$(u_n)_n$ converges to a function $u_0$ in $H^1_{0,r}(\Omega)$;
\item [$b)$]
$|\beta(u_n)|\lo +\infty$.
\end{itemize}

\end{prop}

\proof
Since $(u_n)_n$ is a Palais-Smale sequence for $E$, in particular it
is bounded in $H_{0,r}^1(\Omega)$, so there exists $u_0\in H_{0,r}^1(\Omega)$
such that, up to a subsequence,
\begin{equation}
\label{PS1}
u_n\lo u_0\quad\mbox{weakly in }H_0^1(\Omega)\mbox{ and in } 
L^p(\Omega),\ \mbox{a.e. in }\Omega \mbox{ and in } L^p_{\loc}(\R^N).
\end{equation}
Furthermore, there exists a sequence $(\mu_n)_n$ in $\R$ such that
\begin{equation}
\label{PS2}
(\D E(u_n),w)=\int_\Omega (\D u_n\cdot\D w+u_n w)=\mu_n\int_\Omega
|u_n|^{p-2}u_n w+o(1)\|w\|\quad\forall w\in H_{0,r}^1(\Omega).
\end{equation}
Setting $w=u_n$ in (\ref{PS2}), we get
\begin{equation}
\label{PS2a}
\lim_{n\to +\infty}\mu_n=c.
\end{equation}
So, $u_0$ verifies
\begin{equation}
\label{PS3}
\int_{\Omega} (\D u_0\cdot \D w+u_0 w)=c\int_\Omega |u_0|^{p-2}u_0 w
\qquad \forall w\in H_{0,r}^1(\Omega)
\end{equation}
and, taking into account the principle of symmetric criticality, it
is a solution of
\begin{equation}
\label{PS4}
\left\{
\begin{array}{l}
-\Delta u+u=c|u|^{p-2}u\quad\mbox{ in }\Omega\\
u\in H_0^1(\Omega).
\end{array}\right.
\end{equation}
Now, set
\begin{equation}
\label{PS5}
v_n(x)=u_n(x)-u_0(x).
\end{equation}
If $v_n\to 0$ in $H_0^1(\Omega)$, we are done, taking into account the
continuity of $\beta$, otherwise there exists a constant $k_0>0$ such
that, up to a subsequence,
\begin{equation}
\label{PS9}
\|v_n\|\ge k_0>0\qquad\forall n\in\N.
\end{equation}
We are proving that in such a case $u_0\equiv 0$ and $(b)$ holds.

The sequence $(v_n)_n$ is uniformly bounded in $H_{0,r}^1(\Omega)$ and, by
(\ref{PS1}),
\begin{equation}
\label{PS6}
v_n\lo 0\quad\mbox{weakly in }H_0^1(\Omega)\mbox{ and in } 
L^p(\Omega),\ \mbox{a.e. in }\Omega \mbox{ and in } L^p_{\loc}(\R^N);
\end{equation}
moreover, a direct computation shows that
\begin{equation}
\label{PS7}
\|v_n\|^2=\|u_n\|^2-\|u_0\|^2+o(1)
\end{equation}
and, by the Brezis-Lieb Lemma (\cite{BLi}),
\begin{equation}
\label{PS8}
|v_n|_p^p=|u_n|_p^p-|u_0|_p^p+o(1).
\end{equation}
Observe that, by (\ref{PS2}),(\ref{PS2a}),(\ref{PS3}) and (\ref{PS7}),
\begin{equation}
\label{PS10}
\|v_n\|^2=c(|u_n|_p^p-|u_0|_p^p)+o(1);
\end{equation}
hence, from (\ref{PS8}), $|v_n|_p^p\ge k_1>0$ follows, for a suitable
constant $k_1>0$. 

Now, let us decompose $\R^N$ into the $N$-dimensional hypercubes
$Q_l$, having unitary sides and vertices with integer co-ordinates and
put, for all $n\in\N$,
\begin{equation}
\label{PS11}
d_n=\max_{l\in\N}|v_n|_{L^p(Q_l)}.
\end{equation}

We claim that there exists $\gamma>0$ such that
\begin{equation}
\label{PS12}
d_n\ge\gamma>0\qquad\forall n\in\N;
\end{equation}
indeed
\begin{eqnarray}
\vspace{1mm}
& 0<k_1& \displaystyle{\le|v_n|_p^p=\sum_{l\in\N} |v_n|_{L^p(Q_l)}^p\le
\max_{l\in\N}|v_n|_{L^p(Q_l)}^{p-2}\sum_{l\in\N}
|v_n|_{L^p(Q_l)}^2}\nonumber \\
& & \displaystyle{
\le
d_n^{p-2} s\sum_{l\in\N} \|v_n\|_{H^1(Q_l)}^2=
d_n^{p-2} s \|v_n\|^2,}\label{PS13}
\end{eqnarray}
where $s$ is given by the embedding $H^1(Q_l)\hookrightarrow L^p(Q_l)$
(it is indipendent of $l$).
 
For all $n\in\N$, let $y_n=(y_n',y_{n,N})$ be the center of an
hypercube $Q_n$ where
\begin{equation}
\label{PS14}
|v_n|_{L^p(Q_n)}\ge \gamma.
\end{equation}
Since the functions $(v_n)_n$ are radially symmetric with respect to
the $x_N$-axis and uniformly bounded in $L^p(\R^N)$, by using (\ref{PS14})
we get 
\begin{equation}
\label{f1}
|y'_n|\le k,
\end{equation}
for a suitable constant $k>0$. Hence we must have
\begin{equation}
\label{PS15}
|y_{n,N}| \lo +\infty,
\end{equation}
because $v_n\to 0$ in $L^p_{\loc}(\R^N)$.

Now, let us set $\tilde y_n=(0,\ldots,0,y_{n,N})$ and call $\tilde
v_0$ the weak limit, in $H^1 (\R^N)$ and in 
$L^p(\R^N)$, of $\tilde v_n(x)=v_n(x+\tilde y_n)$, up to a subsequence.
From (\ref{PS14}), (\ref{f1})  and Rellich Theorem we infer $\tilde v_0\not\equiv
0$ and, taking into account (\ref{PS2}) and (\ref{PS2a}), $\tilde v_0$
solves problem $P(c,\cD)$, on its domain $\cD$ (see(\ref{PmD})).

So, it is possible to conclude that $u_0\equiv 0$.
Indeed (\ref{PS7}) implies $\|u_n\|^2\ge \|u_0\|^2+\|\tilde
v_0\|^2+o(1)$, thus, if $u_0$ were not zero, by Proposition \ref{C1}
applied to $u_0$ and $\tilde v_0$ we would obtain
\[
E(u_n)=\|u_n\|^2>2m\left({m\over c}\right)^{2\over p-2}+o(1),
\]
from which the limit value $c$ verifies $c>2^{1-{2\over p}}m$,
contrary to our assumption.

Now, we claim that
\begin{equation}
\label{n1}
\tilde v_n\lo\tilde v_0\qquad\mbox{ in } H_0^1(\cD).
\end{equation}
In fact, if this is not the case, we can argue for $(\tilde v_n)_n$ as
we have done for the sequence $(u_n)_n$, obtaining that $c\ge
2^{1-{2\over p}}m$, a contradiction. 
From (\ref{n1}), and taking into account (\ref{PS15}), we infer that
$(b)$ holds. 

\qed

\begin{rem}
\label{R}
{\em
Actually, Proposition \ref{PS} holds if $c=2^{1-{2\over p}}m$, too.
In fact, in such a case, first one proceed as in the proof of
Proposition \ref{PS} and then a further analysis of the Palais-Smale
sequence $(u_n)_n$ shows that $c=2^{1-{2\over p}}m$ is possible only if
the limit function $u_0$ is identically zero and $(u_n)_n$ tends to
split in two different waves, infinitely distant one from the other and such
that each one approaches the function $\omega$ on a limit domain that
has to be $\R^N$.
So one can conclude that if $\beta(u_n)=(0,\ldots,0,\zeta_n)$ then $\zeta_n\lo +\infty$ must hold, by the
shape of $\Omega$ in the half space $\{x_N<0\}$.

This means that when $c=2^{1-{2\over p}}m$ a Palais-Smale sequence is
  composed by two waves ``almost minimizing'' on strips infinitely
  distant each other and whose sizes become infinite.
}
\end{rem}


\sezione{Proof of the main result}


\begin{lemma}
\label{L1}
Set
\begin{equation}
\label{bi}
\cB:=\inf\left\{E(u)\ :\ \beta(u)={Q^+_i+Q_{i+1}^-\over 2}\mbox{ or
}\beta(u)=-{Q^+_i+Q_{i+1}^-\over 2}\mbox{ for some }i\in\N\right\},
\end{equation}
it holds
\[
\cB >m.
\]
\end{lemma}
\proof
Let us define
\[
\widetilde\cD=\{x=(x',x_N)\in\R^N\ :\ |x'|<1\}\cup\left\{x=(x',x_N)\ :\
|x_N|>{1\over 2}\right\}.
\]
Observe that $\Omega\subset \widetilde\cD+{Q^+_i+Q_{i+1}^-\over 2}$
and  $\Omega\subset \widetilde\cD-{Q^+_i+Q_{i+1}^-\over 2}$, $\forall
i\in\N$, hence, taking into account (\ref{b1}), if we prove that 
\begin{equation}
\label{G}
\inf\{\|u\|^2\ :\ u\in H_{0,r}^1(\widetilde\cD),\ |u|_p=1,\ \beta(u)=0\}>m
\end{equation}
we are done.

Suppose, by contradiction, that (\ref{G}) does not hold. 
Then there exists a sequence $(u_n)_n$ in $H_{0,r}^1(\widetilde\cD)$,
$|u_n|_p=1$ $\forall n\in\N$, such that
\begin{equation}
\label{el1}
\lim_{n\to +\infty}\|u_n\|^2=m,
\end{equation}
\begin{equation}
\label{el2}
\beta(u_n)=0\qquad\forall n\in\N.
\end{equation}

Since the sequence $(u_n)_n$ is in $H_{0,r}^1(\widetilde\cD)$ and is
minimizing for (\ref{em}), then (see \cite{BC}, for example)  there
exist a sequence of points $(y_n)_n$ in $A$ and a 
sequence of axially symmetric functions $(w_n)_n$ in $H^1(\R^N)$ such that
\begin{equation}
\label{el3}
u_n(x)=\omega (x-y_n)+w_n(x)\qquad\mbox{ with }\lim_{n\to +\infty}
\|w_n\|=0.
\end{equation}
From (\ref{b1}), (\ref{b2}) and the continuity of $\beta$ it follows
that
\begin{equation}
\label{el4}
\lim_{n\to +\infty}|\beta(u_n)-y_n|=
\lim_{n\to +\infty}|\beta(\omega (x-y_n)+w_n(x))-y_n|=
\lim_{n\to +\infty}|\beta(\omega (x)+w_n(x+y_n))|=0.
\end{equation}
By (\ref{el2}), (\ref{el4}) and (\ref{el3}) we have that
$u_n\to \omega$ in $H^1(\R^N)$; this is not possible, because
$u_n = 0$ on $\R^N\setminus\widetilde\cD$, $\forall n\in \N$, and $\omega>0$ on
$\R^N$, so the proof is completed.

\qed

\begin{lemma}
\label{D}
Let $k_1\leq k_2<\min\{\cB,2^{1-{2\over p}}m\}$ and assume that
\begin{equation}
\label{D1}
\{u\in V\ :\ E(u)\in [k_1,k_2],\ \D E(u)=0\}=\emptyset;
\end{equation}
if $\bar u\in E^{k_2}$ then there exist $\e>0$ and a continuous path
$\eta:[0,1]\to E^{k_2}$ such that
\begin{equation}
\label{D2}
\eta(0)=\bar u,\quad \eta(1)\in  E^{k_1-\e},\quad t\mapsto
E(\eta(t))\mbox{ is decreasing,}
\end{equation}
where  $E^b$, $b\in\R$, denotes the sublevel 
\[
E^b=\{u\in V\ :\ E(u)\le b\}.
\]
\end{lemma}

\proof
By assumption (\ref{D1}), $\bar u$ has a neighborhood $\cU(\bar u)$ in $V$
such that
\[
\|\D E (u)\|>\const>0\qquad\forall u\in\cU(\bar u).
\]
Then, the Cauchy problem
\[
(C) \hspace{1cm}
\left\{
        \begin{array}{l}
        {d\phantom{i}\over dt}\, \Phi (t) =-{\D E(\Phi(t))\over \|\D E(\Phi(t))\|^2} \\
 \Phi(0)=\bar u
        \end{array}
\right.
\]
can be locally solved, on $V$.
Let $[0,T)$ be the maximal interval where $\Phi$ can be defined; we shall
  show that $L:=\lim\limits_{t\to T^-} E(\Phi(t))$ is well defined and
  $L< k_1$.

First, we claim that $t\mapsto E(\Phi(t))$ is a decreasing function;
indeed, for $t_1,t_2\in [0,T)$ we have
\begin{equation}
\label{D3}
E(\Phi(t_2))-E(\Phi(t_1))=\int_{t_1}^{t_2}{d\phantom{i}\over dt}\,
E(\Phi (t))\,dt
=\int_{t_1}^{t_2}\D E(\Phi (t))\cdot {d\phantom{i}\over dt}\, \Phi
(t)\, dt=-(t_2-t_1).
\end{equation}
From (\ref{D3}), and $(a)$ of Proposition \ref{MM}, it follows, in
particular, that $T<+\infty$.

Now, let us suppose, by contradiction, that $L\ge k_1$ and let
$(t_n)_n$ be a sequence in $[0,T)$ such that $t_n\to T$
  and $\|\D E(\Phi(t_n))\|\lo 0$, as $n\to +\infty$.
  Such a sequence $(t_n)_n$ must exist, otherwise $\Phi$ would be
  Lipshitz, because it solves $(C)$, and so we could extend $\Phi$ on an
  interval $[0,T')$ with $T'>T$, contrary to the maximality of $[0,T)$.

By applying Proposition \ref{PS} to $\Phi(t_n)$, up to a subsequence,
we have two cases: 

\vspace{1mm}

\hspace{3mm} $\Phi(t_n)$ converges to a function $\bar v\in V$; then
\[
E(\Phi(t_n))\lo E(\bar v)=L\in [k_1,k_2],\qquad \D E(\Phi(t_n))\lo  \D
E(\bar v) =0,
\]
contrary to (\ref{D1}).

\vspace{1mm}

\hspace{3mm} $|\beta(\Phi(t_n))|\lo +\infty$; then, since $\beta$ is
  continuous, there exists $\bar t\in [0,T)$
  such that $\beta(\Phi(\bar t))={Q_i^++Q_{i+1}^-\over 2}$ or
  $\beta(\Phi(\bar t))=-{Q_i^++Q_{i+1}^-\over 2}$, for some 
  $i\in \N$,  and, moreover,
  $E(\Phi(\bar t))\le E(\Phi(0))<\cB$, by assumption. 
This is not possible by definition of $\cB$ (see (\ref{bi})), so we
  have the desired conclusion, up to a normalization of the path.

\qed

{\noindent{\underline {\sf Proof of Theorem \ref{T}} \hspace{2mm}}}
Observe that, taking into account Corollary \ref{C2}, the critical points of $E$ whose
critical levels are in $(m,2^{1-{2\over p}}m]$ are
nonnegative functions, actually positive by the maximum principle.

Then, in order to prove Theorem \ref{T}, it is sufficient to show that for
every $a\in(m,2^{1-{2\over p}}m]$ there exists a critical value $c$ for
$E$ such $c\in (m,a)$.

Let $S_k$ be a ``strip'' such that $\|z_k\|^2<\min\{a,\cB\}$
(see (\ref{S}), (\ref{q1}) and (\ref{q00})); we claim that $E$ has a
critical value $c\in (m,\|z_k\|^2]$. 
Assume, by contradiction, that no critical value exists in $[m,
 \|z_k\|^2]$. 
Then, by Lemma \ref{D}, there exists a continuous path from
$z_k\left(x-{Q_k^++Q_k^-\over 2}\right)$ to a function $v\in V$ such
that $E(v)< m$. But this contradicts (\ref{eC5}), so we get the claim,
by Remark \ref{R5}.  

\qed



\end{document}